\documentclass[11pt]{article}
\usepackage{latexsym,amssymb}

\def\demo{\noindent{\bf Proof. }}
\def\QED{\hfill$\Box$}
\def\qed{\QED}

\newtheorem{Theorem}{Theorem}[section]
\newtheorem{Lemma}[Theorem]{Lemma}
\newtheorem{Corollary}[Theorem]{Corollary}
\newtheorem{Proposition}[Theorem]{Proposition}
\newtheorem{Remark}[Theorem]{Remark}
\newtheorem{Example}[Theorem]{Example}
\newtheorem{Conjecture}[Theorem]{Conjecture}
\newtheorem{Definition}[Theorem]{Definition}

\def\rk{\mbox{\rm rank}}
\def\pp{\mathbb P}
\def\zz{\mathbb Z}
\def\qq{\mathbb Q}
\def\rar{\rightarrow}
\def\lar{\longrightarrow}

\begin{document}
\topmargin3mm
\hoffset=-2.3cm
\voffset=-1.5cm

\

\vskip 3cm


\begin{center}

{\large\bf
Linear syzygies and birational combinatorics}\vspace{6mm}\\
{\normalsize Aron Simis and Rafael H. Villarreal*}
\end{center}

\date{}

\medskip

\begin{abstract}
\noindent
Let $F$ be a finite set of monomials of the same
degree $d\geq 2$ in a polynomial ring $R=k[x_1,\ldots,x_n]$ over an
arbitrary field $k$. We give some
necessary
and/or sufficient conditions for the birationality of the ring
extension
$k[F]\subset R^{(d)}$, where $R^{(d)}$ is the $d${\it th} Veronese
subring
of $R$.
One of our results extends to arbitrary characteristic, in the case of
rational
monomial maps, a previous syzygy-theoretic birationality criterion
in characteristic zero obtained in
\cite{CRS}.
\end{abstract}

\section{Introduction}

By the expression ``birational combinatorics'' we mean the theory of
characteristic-free rational maps $\pp^{n-1}\dasharrow \pp^{m-1}$ defined by
monomials, along with natural criteria for such maps to be birational onto their
image varieties. Both the theory and the criteria are intended to be simple and
typically reflect the monomial data, as otherwise one falls back in the general
theory of birational maps in projective spaces (cf., e.g., \cite{RuSi},
\cite{Si2}).

A first incursion in this kind of theory was made in \cite{SiVi}. There one focused
mainly on monomial rational maps whose base ideal (ideal theoretic base locus) was
normal. Though the results were fairly complete and some of the techniques used
there are repeated here, one felt that normality was a special case obscuring the
general picture.

In the present paper we envisage a general theory focusing on the underlying
combinatorial elements rather than on special algebraic properties of the base
ideal. In this sense, what we accomplish goes in the opposite direction of recent
work on birational maps, where the emphasis fell on special behavior of the base
locus. On the other hand, we did draw upon \cite{RuSi} and \cite{Si2} (also upon
the ongoing \cite{CRS}) by invoking the role played by the so-called linear
syzygies of the coordinates of the rational map. The methods in the first two of
these references are specially suited for the explicit computation of the inverse
map of a birational map onto the image. To compromise between the two approaches,
we show a bridge between them by means of comparing the respective linear algebra
gadgets - from modules over the ground polynomial ring $k[x_1,\ldots,x_n]$ to
modules over $\zz$. The challenge remains as to how one computes the inverse map by
a purely combinatorial method.

We now describe the content of the paper in more detail. It goes without saying
that the language throughout is algebraic or combinatorial, although we do add
frequent remarks as to the geometric meaning of the results.

Section~\ref{olinda} sets up the scenario for the basic pertinent integer
combinatorics. We emphasize two criteria of birationality in this setup - the {\it
arithmetical principle of birationality\/} and the {\it determinantal principle of
birationality}. These criteria were used in \cite{SiVi} and seem to be part of the
folklore in the scattered literature. Then, we introduce the various versions of
matrices that will play a distinctive role in the theory and, in particular, replay
in more generality the passage from the transposed Jacobian matrix to the
log-matrix of a set of monomials, as devised in \cite{Si1}. Since we wish to remain
characteristic-free, we take the formal Jacobian matrix rather than the ordinary
one, as is explained in the section. The so-obtained numerical matrices allow for a
first birationality criterion (Proposition~\ref{ufpe}). We then proceed to a full
arithmetical characterization of birationality (Theorem~\ref{maincriterion}).

Section~\ref{fitting} deals with the role of the Fitting ideals of monomial
structures. We expand on the topic only enough in order to compare ranks between
matrices over $k[x_1,\ldots,x_n]$ and matrices over $\zz$. As a side bonus, we
characterize totally unimodular log-matrices in terms of Fitting ideals of the
formal Jacobian matrix. The main result of the section is Theorem~\ref{boa-viagem},
which extends one of the results of \cite{CRS} to all characteristics for monomial
rational maps.

Section~\ref{degree2} focuses on the case of monomials of degree $2$. Here, we give
complete results, covering all previously known results and establishing facts that
do not extend to higher degrees. We introduce the notion of {\it cohesiveness\/}
for rational maps of any degre inspired by the graph theoretic concept of
connectedness. We show, preliminarily, that the lack of cohesiveness is an early
obstruction for birationality and for the existence of ``enough'' linear syzygies.
If, moreover, the degree is $2$ we show that cohesiveness is a necessary and
sufficient condition for having a linear syzygy matrix of maximal rank
(Proposition~\ref{cohesive_vs_linear}). We proceed to one of the main theorems of
the section (Theorem~\ref{aronvila05}) saying that a rational map of degree $2$ is
birational onto its image if and only if it is cohesive and the corresponding
log-matrix has maximal rank. This comes to us as a bit of a surprise as it says
that any cohesive coordinate projection of the $2$-Veronesean that preserves
dimension is birational onto the image; moreover, this holds in any characteristic.
We have not met any explicit mention of this fact in the previous literature. We
give examples to show how easily this fails for non-monomial rational maps and for
monomial ones in degrees $\geq 3$. Finally we care to translate the results into
the language of graphs with loops.

The last section has the purpose of describing sufficiently ample classes of
monomial rational maps that are birational. It is further subdivided in two
subsections, the first of which is entirely devoted to classes of Cremona maps. We
characterize Cremona transformations of degree $2$ as those cohesive ones whose
log-determinant is nonzero. The corresponding graph theoretic characterization is
suited to construct other Cremona transformations of higher degree via a certain
duality principle. The second subsection is a pointer to a recently studied class
of combinatorial objects called polymatroidal monomial sets. This class includes
the toric algebras of Veronese type which, from the geometric angle, constitutes a
vast class of dimension preserving projections of the ordinary Veronese embeddings.

\section{Birationality of monomial subrings}\label{olinda}

Let $R=k[\mathbf{x}]=k[x_1,\ldots,x_n]$ be a polynomial ring over a field $k$. As
usual we set ${x}^{\alpha}:= x_1^{a_1}\cdots x_{n}^{a_{n}}$ if $\alpha=(a_1,
\ldots,a_n)\in {\mathbb N}^n$. In the sequel we consider a finite set of distinct
monomials ${F}=\{{x}^{v_1},\ldots, {x}^{v_q}\}\subset R$ of the same degree $d\geq
2$ and having no non-trivial common factor. We also assume throughout that $F$ is
not {\it conic}, i.e., that every $x_i$ divides at least one member of $F$. By
trivially contracting to less variables, any set of monomials can be brought to
this form.

Two integer matrices naturally associated to $F$ are:
$$
A=(v_1,\ldots,v_q) \ \mbox{ and }\ A'=
\left(\hspace{-2mm}
\begin{array}{ccc}
v_1&\cdots&v_q\\
1&\cdots&1
\end{array}\hspace{-2mm}
\right),
$$
where the $v_i$'s are regarded as column vectors. We will often refer to $A$ as the
{\it log-matrix\/} of ${F}$.

\medskip

If $C$ is an integer matrix with $r$ rows, we denote by $\mathbb{Z}C$ (resp.
$\mathbb{Q}C$) the subgroup of $\mathbb{Z}^r$ (resp. subspace of $\mathbb{Q}^r$)
generated by the columns of $C$. $\Delta_r(C)$ will denote the greatest common
divisor of all the nonzero $r\times r$ minors of $C$.

An extension $D'\subset D$ of integral domains is said to be birational if it is an
equality at the level of the respective fields of fractions. In the sequel let
${\bf x}_d$ denote the set of {\it all\/} monomials of degree $d$ in $R$. Then
$k[{\bf x}_d]$ is the $d${\it th} Veronese subring $R^{(d)}$ of $R$. Our main aim
is the birationality of the ring extension $K[F]\subset k[{\bf x}_d]$.

\medskip

For convenience of reference, we quote the following easy results stated in
\cite{SiVi}:

\begin{Lemma}\label{principle}{\rm (Arithmetical Principle of Birationality (APB))}
Let $F$ and $G$ be finite sets of monomials of $R$ such that $F\subset G$, and let
$A,B$ be their respective log-matrices. Then $k[F]\subset k[G]$ is a birational
extension if and only if $\mathbb{Z}A= \mathbb{Z}B$.
\end{Lemma}
\demo In this situation, the ring extension is birational if and only every
monomial of ${G}$ can be written as a fraction whose terms are suitable power
products of the monomials of ${F}$. Clearing denominators of such a fraction and
taking log of both members establishes the required equivalence. \qed

\begin{Lemma}\label{principle2}{\rm (Determinantal Principle of Birationality (DPB))}
Let $F$ be a finite set of monomials of the same degree $d\geq 1$. Then
$k[{F}]\subset k[{\bf x}_d]$ is a birational extension if and only if
$\Delta_n(A)=d$.
\end{Lemma}
\demo See \cite[Proposition 1.2]{SiVi}. \qed

\bigskip

By $e_i, 1\leq i\leq n$ we denote the canonical basis vectors of the vector space
${\mathbb R}^n$ (sometimes of the free module $\zz^n$, respectively, the
$\qq$-vector space $\qq^n$). Let, as before,  ${F}=\{{x}^{v_1},\ldots,
{x}^{v_q}\}\subset R$ be a set of monomials of the same degree $d\geq 2$.

\medskip

Consider the following basic matrices:
\begin{description}
\item{\rm (a)} the matrix ${\cal L}{\cal S}(F)$ of the so-called {\it linear
syzygies\/} of $F$, whose columns are the set of vectors of the form
$x_ie'_j-x_ke'_\ell$ such that $x_ix^{v_j}=x_kx^{v_\ell}$;
\item{\rm (b)} the {\it numerical linear
syzygy matrix\/} $S$ obtained from  ${\cal L}{\cal S}(F)$ by making the
substitution $x_i=1$ for all $i$;
\item{(c)} the matrix $M$ whose columns are the set of
difference vectors $e_i-e_k$ such that $e_i-e_k=v_j-v_\ell$, for some pair of
indices $j,\ell\in \{1,\ldots,q\}$ -- in other words,  $M=AS$;
\item{\rm (d)} the {\it formal\/} Jacobian matrix $$\Theta(F)=
\left(\frac{\partial x^{v_j}}{\partial x_i}\right)_{\begin{array}{c} 1\leq j\leq q\\
1\leq i\leq n
\end{array}}$$
\end{description}
A word in order to explain the last matrix. The notion of derivative of a
polynomial $f\in k[x_1,\ldots,x_n]$ usually requires the specification of a base
field. However, if $f$ is an ordinary monomial ${x}^{a}=x_1^{a_1}\cdots x_n^{a_n}$
its {\it formal\/} partial derivative with respect to $x_i$ is defined to be
$$
a_ix_1^{a_1}\cdots x_i^{a_{i-1}} x_i^{a_i-1}x_i^{a_{i+1}}\cdots
x_n^{a_n}
$$
regarded as a term in the polynomial ring ${\mathbb Z}[{\bf x}]$ (in particular it
is always nonzero provided $a_i\geq 1$). The {\it formal Jacobian matrix\/} of
$\,{x}^{v_1},\ldots,{x}^{v_q}$ is accordingly defined. Of course, by applying the
unique homomorphism from $\zz$ to $k$ we find the ordinary partial derivatives and
the ordinary Jacobian matrix over this ring.

\medskip

Notice that the matrices in the first row of the diagram:
$$
\begin{array}{ccccc}
{\Theta}(F)^t&\ \ \ &{\cal L}{\cal S}(F)&\ \ \ &{\cal
M}:={\Theta}(F)^t{\cal L}{\cal S}(F)\\
\downarrow& &\downarrow &
&\downarrow\\
A& & S & &M:=AS \\
\end{array}
$$
specialize to the matrices in the second row by making $x_i=1$ for all $i$. The
matrices ${\cal L}{\cal S}(F)$ and $S$ have order $q\times r$, while the matrices
${\cal M}$ and $M$ have order $n\times r$.

\medskip

Here is a couple of uses of these matrices. The following notion will be used in
the proof below:  a matrix $C$ is called {\it totally unimodular\/} if each
$i\times i$ minor of $C$ is $0$ or $\pm 1$ for all $i\geq 1$.

\begin{Proposition}\label{ufpe} Let $F$ be a finite set of monomials of the
same degree $d\geq 2$.
\begin{enumerate}
\item[{\rm (i)\ }] If\, ${\rm rank}(M)=n-1$, then $k[F]\subset k[{\bf x}_d]$
is a birational extension.
\item[{\rm (ii)}] If
${\rm rank}(S)=q-1$ and ${\rm rank}(A)=n$, then ${\rm rank}(M)=n-1$. In particular
$k[{F}]\subset k[{\bf x}_d]$ is birational.
\end{enumerate}
\end{Proposition}
\demo (i) Let $a=(a_i)\in\mathbb{N}^n$ such that $|a|=\sum_ia_i=d$. By APB
(Lemma~\ref{principle}) it suffices to prove that $a\in\mathbb{Z}A$. Let
$w_1,\ldots,w_r$ be the column vectors of the matrix $M$. Each $w_m$ is of the form
$e_i-e_k=v_j-v_\ell$ for a unique pair $i\neq k, \,1\leq i<k\leq n$ and suitable
$j\neq \ell, \,1\leq j<\ell\leq q$. Hence ${\rm rank}(A)=n$ because
$v_1\notin\mathbb{Q}M$. Therefore we can write
$$
{\lambda}a={\lambda}_1w_1+\cdots+{\lambda}_rw_r+\mu v_1\ \ \
({\lambda},\mu,{\lambda}_i\in\mathbb{Z}).
$$
Taking inner product with $\mathbf{1}=(1,\ldots,1)$ yields
\begin{eqnarray}
\lefteqn{ {\lambda}d={\lambda}|a|={\lambda}_1|w_1|+\cdots+{\lambda}_r|w_r|+
\mu|v_1|=\mu d\ \Rightarrow\ {\lambda}=\mu}&\ \ \ \
 \ \ \ \ \ \ \ \ \ \ \ \ \ \ \ \ \ \ \ \
\ \ \ \
\ &\nonumber\\
&&\Rightarrow \ {\lambda}(a-v_1)={\lambda}_1w_1+\cdots+{\lambda}_rw_r.
\label{ufpe1}
\end{eqnarray}
Consider the digraph $\cal D$ with vertex set $X=\{x_1,\ldots,x_n\}$ such that the
directed edges $(x_i,x_k)$ correspond bijectively to the column vectors $e_i-e_k$
of $M$. The incidence matrix of $\cal D$ is $M$, thus $M$ is totally unimodular
\cite[p.~274]{Schr} and $\mathbb{Z}^n/\mathbb{Z}M$ is torsion-free. Hence from
Eq.~(\ref{ufpe1}) we get $a-v_1\in\mathbb{Z}M$ and
$a\in\mathbb{Z}M+v_1\subset\mathbb{Z}A$, as required.

(ii) Consider the $\qq$-linear maps
$$
\mathbb{Q}^r\stackrel{S}{\longrightarrow}\mathbb{Q}^q
\stackrel{A}{\longrightarrow}\mathbb{Q}^n.
$$
Letting $A_1$ denote the restriction of $A$ to ${\rm im}(S)$, we have a linear map
$$
{\rm im}(S)\stackrel{A_1}\longrightarrow {\rm im}(AS)={\rm im}(M)\longrightarrow 0.
$$
By hypothesis, $\dim({\rm im}(A))=n$ and $\dim({\rm im}(S))=q-1$. Hence
\begin{eqnarray*}
q-1&=&\dim({\rm im}(S))=\dim({\rm ker}(A_1))+{\rm dim}({\rm
im}(M)),\\
q-n&=&{\rm dim}({\rm ker}(A))\geq{\rm dim}({\rm ker}(A_1)).
\end{eqnarray*}
Therefore ${\rm dim}({\rm im}(M))\geq n-1$. On the other hand, since ${\rm im}(M)$
is generated by vectors of the form $e_i-e_j$, certainly $e_1\not\in {\rm im}(M)$,
hence ${\rm dim}({\rm im}(M))=n-1$. \QED

\begin{Remark}\label{apr14-05}\rm Let $\cal D$ be the digraph in the
proof of Proposition~\ref{ufpe}(i). Then according to
\cite[Theorem~8.3.1]{godsil} we have
$$
{\rm rank}(M)=n-c,
$$
where $c$ is the number of connected components of $\cal D$. In
particular $M$ has rank $n-1$ if and only if $\cal D$ is connected.
\end{Remark}

To proceed with a full arithmetical characterization of birationality, we will need
the following results on modules over $\zz$.

\begin{Lemma}\label{lemmata} {\rm (i)} Let $e_i,1\leq j \leq n$ be the canonical basis
vectors of the free $\zz$-module $\zz^n$ and let $E\subset \zz^n$ be the submodule
generated by the difference vectors $e_i-e_k, 1\leq i< k \leq n$. Then $E$ is
freely generated by  $\{e_1-e_k\,|\, 2\leq k \leq n\}$ and the quotient $\zz^n/E$
is torsionfree of rank one.
\newline {\rm (ii)} Let  $\alpha_1,\ldots,\alpha_m\in \zz^n$ be arbitrarily given.
Then the injective $\zz$-homomorphism $\zz^n\rar\zz^{n+1}$, $\alpha\mapsto
(\alpha,0)$, induces an injective homomorphism of $\zz$-modules
$$
\mathbb{Z}^n/\zz\left(\alpha_2-\alpha_1,
\ldots,\alpha_m-\alpha_1\right)\hookrightarrow
\mathbb{Z}^{n+1}/\zz\left((\alpha_1,1), \ldots,(\alpha_m,1)\right),
$$
which is an isomorphism at the level of the respective torsion submodules.
\end{Lemma}
\demo (i) This is simply the fact that $E$ is the kernel of the $\zz$-homomorphism
$\zz^n\rar\zz$, $(a_1,\ldots,a_n)\mapsto a_1+\cdots +a_n$.

(ii) Clearly, there is an induced map as argued -- because $\alpha_j-\alpha_1$ maps
to $(\alpha_j-\alpha_1,0)=(\alpha_j,1)-(\alpha_1,1)$ -- and the induced map is
injective -- because the two equations $a_1+\cdots +a_m=0$ and
$\alpha=a_1\alpha_1+\cdots +a_m\alpha_m$ easily imply that $\alpha\in
\zz\left(\alpha_2-\alpha_1, \ldots,\alpha_m-\alpha_1\right)$.

Next, clearly any homomorphism maps torsion to torsion, so it remains to check
surjectivity at the torsion level.  Let then $(\alpha,b)$ be a torsion element of
$\mathbb{Z}^{n+1}/\zz((\alpha_1,1), \ldots,(\alpha_m,1))$. This implies a relation
$$
s(\alpha,b)=\lambda_1(\alpha_1,1)+\cdots+ \lambda_m(\alpha_m,1)\ \ \ \ \ \
(\lambda_i\in\mathbb{Z}),
$$
where $0\neq s\in\mathbb{N}$, $\alpha\in\mathbb{Z}^n$ and $b\in\mathbb{Z}$. Then
\begin{eqnarray*}
s\alpha&=&\lambda_1\alpha_1+\cdots+\lambda_m\alpha_m,\\
sb&=&\lambda_1+\cdots+\lambda_m,\\
s(\alpha-b\alpha_1)&=& \lambda_2(\alpha_2-\alpha_1)+\cdots+
\lambda_m(\alpha_m-\alpha_1).
\end{eqnarray*}
Hence it follows that the class $\overline{\alpha-b\alpha_1}\in
\zz^n/\zz(\alpha_2-\alpha_1, \ldots,\alpha_m-\alpha_1)$ is a torsion element and
maps to $\overline{(\alpha,b)}$, as required. \QED

\begin{Theorem}\label{maincriterion} Let $F$ be a finite set of monomials of the
same degree $d\geq 2$. The following conditions are equivalent
\begin{description}
\item{\rm (a)} $k[{F}]\subset k[{\bf x}_d]$ is
birational.
\item{\rm (b)} $\mathbb{Z}^n/\mathbb{Z}(\{v_1-v_j\vert\,
2\leq j\leq q\})$ is free of rank $1$.
\item{\rm (c)} The log-matrix $A$ of $F$  has maximal rank and
$\mathbb{Z}(\{v_1-v_j\vert\, 2\leq j\leq q\})=\mathbb{Z}(\{e_1-e_k\vert\, 2\leq
i\leq n\})$.
\end{description}
\end{Theorem}
\demo First we observe that, quite generally, there is an exact sequence of finite
abelian groups
\begin{equation}\label{ufpe2}
0\rar T(\mathbb{Z}^{n+1}/\mathbb{Z}{A}') \stackrel{\varphi}{\lar}
T(\mathbb{Z}^n/\mathbb{Z}A) \stackrel{\psi}{\lar} \mathbb{Z}_d\rar 0
\end{equation}
(here $\varphi((\overline{\alpha,b)})=\overline{\alpha}$ and
$\psi(\overline{\alpha})=\overline{\langle \alpha,{\bf 1}\rangle}$, for
$\alpha\in\mathbb{Z}^n, b\in\mathbb{Z})$ -- see \cite[Proof of Theorem 1.1]{SiVi}.

If, moreover, $A$ has full rank then $\mathbb{Z}^n/\mathbb{Z}A$ is torsion, hence
$\mathbb{Z}^n/\mathbb{Z}A\simeq \zz_d$ if and only if
$\mathbb{Z}^{n+1}/\mathbb{Z}{A}'$ is torsionfree, and in this case
the $0${\it th\/} Fiting
ideal $\Delta_n(A)$ of $\mathbb{Z}^n/\mathbb{Z}A$ is the same as that of $\zz_d$,
i.e., $\Delta_n(A)=(d)$.

On the other hand, we have an exact sequence of $\zz$-modules
\begin{equation}\label{ufpe3}
0\rar \zz A/\zz(\{v_1-v_j\vert\, 2\leq j\leq q\})\rar \zz^n/\zz(\{v_1-v_j\vert
\,2\leq j\leq q\})\rar \zz^n/\zz A\rar 0
\end{equation}
Again, if $A$ has full rank then the leftmost module has rank $1$ and, since the
rightmost module is torsion, the mid module has rank $1$. Now apply
Lemma~\ref{lemmata}(ii) with $m=q$ and $\alpha_j=v_j$ to get
$$T(\zz^n/\zz(\{v_1-v_j\vert \,2\leq j\leq q\}))\simeq T(\zz^{n+1}/\zz A').$$
Therefore, the equivalence (a) $\Longleftrightarrow$ (b) follows from DBP of
Lemma~\ref{principle2}.

It remains to show that (b) $\Longleftrightarrow$ (c). First, (c) $\Rightarrow$ (b)
is clear by Lemma~\ref{lemmata}(i). For the reverse implication, since the mid term
of the sequence (\ref{ufpe3}) is assumed to be torsionfree of rank one and
$\zz(\{v_1-v_j\vert \,2\leq j\leq q\})\neq \zz A$, then $A$ must have full rank
and, moreover, $\zz A/\zz(\{v_1-v_j\vert \,2\leq j\leq q\})$ is torsionfree of rank
one. In particular, there is a splitting $\zz A\simeq \zz(\{v_1-v_j\vert\, 2\leq
j\leq q\})\oplus \zz$ which, after extending to $\qq$, implies
\begin{equation}\label{march28-05}
\mathbb{Q}(\{v_i-v_j\vert 1\leq i<j\leq q\})=\mathbb{Q}(\{e_i-e_j\vert 1\leq
i<j\leq n\}).
\end{equation}
Hence we get the desired equality because of the torsion freeness
hypothesis. Notice that Eq.~(\ref{march28-05}) also follows directly.
Indeed if $\{v_1,\ldots,v_{n}\}$ is a basis for the column space of
$A$, then
$$
\{v_1-v_n,v_2-v_n,\ldots,v_{n-1}-v_n,v_n\}
$$
is also a basis because $|v_i|=d$ for all $i$. Hence each $e_i-e_j$
can be written as
$$
e_i-e_j=a_1(v_1-v_n)+\cdots+a_{n-1}(v_{n-1}-v_n)+a_nv_n\ \ \
(a_i\in\mathbb{Q}).
$$
Taking inner products with the vector $\mathbf{1}=(1,\ldots,1)$ yields
$a_{n}=0$. Therefore we have shown the containment ``$\supset$'' in
Eq.~(\ref{march28-05}). A symmetric argument proves the equality. \QED

\medskip

\section{When are the Fitting ideals monomial ideals?}\label{fitting}

In  \cite[Lemma 1.1]{Si1} was shown that the minors of the Jacobian matrix of a set
of monomials are always monomials (possibly zero). The following result extends and
clarifies the above assertion.

\begin{Proposition}\label{graded_nonsense}
Let $R$ be a graded ring with grading given by an additive abelian monoid ${\cal
Z}$. Let $N$ be a finitely generated ${\cal Z}$-graded module over $R$. Then the
Fitting ideals of $N$ are homogenous ideals of $R$.
\end{Proposition}
\demo By assumption, there is an exact sequence of ${\cal Z}$-graded modules over
$R$
$$\sum_{{\mathfrak z}_j\in {\cal Z}}\, R({\mathfrak z}_j)\stackrel{\phi}{\lar}
\sum_{{\mathfrak w}_i\in {\cal Z}}\, R({\mathfrak w}_i)\lar N\rar 0.
$$
A Fitting ideal of $N$ is an ideal $I_t(\phi)$ generated by the $t$-minors of
$\phi$, for a suitable $t$. This ideal is the image of the well-known induced
${\cal Z}$-graded homomorphism
$$\bigwedge^t \,\sum_{{\mathfrak z}_j\in {\cal Z}}\, R({\mathfrak z}_j)\otimes_R
\bigwedge^t\,\sum_{{\mathfrak w}_i\in {\cal Z}}\, R({\mathfrak w}_i)\lar R.$$
Therefore, $I_t(\phi)$ is a homogeneous ideal of $R$. \qed

\begin{Corollary}\label{fitting_of_monomials}
Let $R=k[x_1,\ldots,x_n]$ be given the standard multigrading {\rm (}i.e., the
$\zz^n$-grading with $x_i$ of degree $(0,\ldots,0,1,0,\ldots,0)${\rm )}. If $N$ is
a finitely generated multigraded $R$-module, then the Fitting ideals of $N$ are
monomial ideals. In particular, any minor of the Jacobian matrix, respectively, of
the syzygy matrix of arbitrary order, of a finite set of monomials is a monomial.
\end{Corollary}
\demo Apply Proposition~\ref{graded_nonsense} while noticing that a homogeneous
polynomial in the standard multigrading is necessarily a monomial. \qed

\medskip

We can also apply the previous result in the case of the standard multigraded ring
$\zz[x_1,\ldots,x_n]$, with $\zz$ in degree ${\bf 0}=(0,\ldots,0)$. The result is
that, in particular, the formal Jacobian matrix of a finite set of monomials has
monomial Fitting ideals. We wish to emphasize this in the following form:

\begin{Corollary}\label{formaljac}
The formal Jacobian matrix and the log-matrix of a finite set $F$ of monomials have
the same number of zero or nonzero minors. In particular, these matrices have the
same rank. Also, there are at most finitely many field characteristics over which
the Jacobian matrix of $F$ over these characteristics has rank strictly smaller
than the rank of the corresponding log-matrix.
\end{Corollary}

There is also a consequence tied up with the notion of a unimodular matrix.

\begin{Corollary} The following are equivalent for a finite set $F$ of monomials.
\begin{enumerate}
\item[{\rm (i)}] The log-matrix of $F$ is totally unimodular
\item[{\rm (ii)}] Every nonzero minor of the formal Jacobian matrix of
$F$ has unit leading coefficient
\item[{\rm (ii)}] The formal Jacobian matrix of $F$
has characteristic-free Fitting ideals {\rm (}i.e., the Fitting ideals of $F$ over
any field are generated by the same set of nonzero monomials{\rm )}.
\end{enumerate}
\end{Corollary}

As for the syzygies of $F$, we observe that, in particular, any minor of the first
Taylor syzygy matrix of $F$ (see \cite{Eisen} for an explanation of
the Taylor complex) is a monomial with coefficient $\pm 1$. We next
include an alternative elementary proof of this fact alone, as the method of the
proof might be useful in some other context.

\begin{Lemma}\label{chicapitanga}
Let ${\cal T}(F)$ denote the Taylor syzygy matrix of $F$. Then any nonzero minor of
${\cal T}(F)$ is a monomial with coefficient $\pm 1$.
\end{Lemma}
\demo We proceed by induction on the size $s$ of the minor. The case $s=1$ being
obvious, we assume that $s\geq 2$. We may clearly assume that the given minor is
formed by the submatrix $Z$ with  the first $s$ rows and columns of ${\cal T}(F)$.
Let $Z'$ denote the $q\times s$ submatrix of ${\cal T}(F)$ with the first $s$
columns. By definition of the Taylor syzygy matrix of $F$, any column of the latter
has exactly two nonzero entries. It follows that the complementary rows in $Z'$ to
the rows of $Z$ cannot all be zero as otherwise $Z$ would be a matrix of syzygies
of the initial $s$ monomials $\{x^{v_1},\ldots,x^{v_s}\}$ of $F$, which is
impossible since $\det(S)\neq 0$ while the entire syzygy matrix of these monomials
has rank $s-1$.

Thus, there
must be a nonzero entry in some complementary row to $Z$ in $Z'$,
say, the
$j$th column,
with $1\leq j\leq s$. By the Taylor construction,  there is exactly
one further
nonzero entry on the $j${\it th} column. This entry must belong to
$Z$ as
otherwise $\det(Z)=0$. Also, this entry is again monomial with
coefficient $\pm$.
Expanding $\det(Z)$
by the $j$th column yields the product of this monomial by the
minor of a suitable
$(s-1)\times (s-1)$ submatrix of $S$. By induction, this minor
has the required
form, hence so does $\det(Z)$.
 \QED

\begin{Corollary}\label{chicapitanga2}
Let $Z$ be any submatrix of ${\cal T}(F)$ and let ${\mathbb T}(Z)$ denote the
specialized matrix over ${\mathbb Z}$ obtained by sending $x_i\mapsto 1$. Then
$\rk(Z)= \rk\, {\mathbb T}(Z)$.
\end{Corollary}

The next result complements one of the
results of \cite{CRS}, where a criterion is
given for a rational map to be
birational in characteristic zero. The present
proposition extends the latter result
in all characteristics for monomial
rational
maps.

\begin{Theorem}\label{boa-viagem} Let $F$ be a finite set of monomials of the
same degree $d\geq 2$. If ${\rm rank}(A)=n$ and ${\rm rank}({\cal L}{\cal
S}(F))=q-1$, then $k[F]\subset  k[{\bf x}_d]$ is birational.
\end{Theorem}
\demo By Corollary~\ref{chicapitanga2} (or by Corollary~\ref{fitting_of_monomials})
the matrix $S$ obtained from ${\cal L}{\cal S}(F)$ by making $x_i=1$ for all $i$
has also rank $q-1$. By Proposition~\ref{ufpe}(ii), the extension $k[F]\subset
k[{\bf x}_d]$ is birational. \QED

\begin{Corollary}\label{apr8-05}
If the log-matrix $A$ has maximal rank and the ideal $I=(F)\subset R$ has a linear
presentation, then $k[F]\subset k[{\bf x}_d]$ is birational.
\end{Corollary}
\demo It follows at once from Theorem~\ref{boa-viagem} because in this case the
rank of ${\cal L}{\cal S}(F)$ is $q-1$. \QED

\section{Monomials of degree two}\label{degree2}

The birational theory of monomials of degree two can be completely established
using elementary graph theory as we show in the sequel.

\medskip

We start with a general auxiliary result which holds, more generally, for any rational
map between projective spaces.

\begin{Lemma}\label{nonsense}
Let $F=\{f_1,\ldots,f_q\}\subset R=k[{\bf x}]=k[x_1,\ldots,x_n]$ be forms of fixed
degree $d\geq 2$. Suppose one has a partition ${\bf x}={\bf y}\cup {\bf z}$ of the
variables such that $F=G\cup H$, where the forms in the set $G$ {\rm
(}respectively, $H${\rm )} involve only the ${\bf y}$-variables {\rm
(}respectively, ${\bf z}$-variables{\rm )}. If neither $G$ nor $H$ is empty then:
\begin{enumerate}
\item[{\rm (i)}] The extension $k[F]\subset k[{\bf x}_d]$ is {\sc not} birational
\item[{\rm (ii)}] The linear syzygy matrix of $F$ does {\sc not} have maximal rank.
\end{enumerate}
\end{Lemma}
\demo  (i) Suppose to the contrary, i.e., that $k(F)= k({\bf x}_d)$. Since clearly
$k(F)=k(G,H)\subset k({\bf y}_d, {\bf z}_d)$, it follows that $k({\bf x}_d)=k({\bf
y}_d, {\bf z}_d)$. Say ${\bf y}=\{y_1,\ldots, y_r\}$ and ${\bf
z}=\{z_1,\ldots,z_s\}$. Then one has $k({\bf y}_d)=k(y_2/y_1,\ldots,y_r/y_1,
y_1^d)$ and, similarly, $k({\bf z}_d)=k(z_2/z_1,\ldots,z_s/z_1, z_1^d)$. But this
is a contradicition as, e.g., $y_1^{d-1}z_1\not\in k(y_2/y_1,z_2/z_1,\ldots,
y_r/y_1,z_s/z_1,y_1^d,z_1^d)$ (for instance, by APB (Lemma~\ref{principle})).

(ii) The linear syzygy matrix of $F$ is a block-diagonal $(r+s)\times m$ matrix
$${\cal L}{\cal S}(F)=\left(\begin{array}{cc}
A&0\\
0&B \end{array} \right), $$ where $A$ and $B$ are the linear syzygy matrices of $G$
and $H$, respectively. Since $\rk(A)\leq r-1$ and $\rk(B)\leq s-1$, then $\rk
({\cal L}{\cal S}(F))\leq r+s-2\leq q-2$. \qed

\begin{Definition}\rm A set $F=\{f_1,\ldots,f_q\}$ of forms of fixed
degree $\geq 2$ will said to be {\it cohesive\/} if the forms have no non-trivial
common factor and $F$ cannot be disconnected as in the hypothesis of the previous
lemma.
\end{Definition}

\begin{Remark}\rm The reason to assume that the forms have no non-trivial
common factor is technical: multiplying a set of forms of the same degree by a
given form yields the same rational map. To make the rational map correspond
uniquely to a set of forms, one usually assumes that their gcd is one, i.e., that
the ideal generated by these forms in the polynomial ring has codimension at least
two (for further details on this and similar matters see \cite{Si2}).
\end{Remark}

Yet another concept that fits the scene is a convenient extension of the notion of
an ideal of linear type.

\begin{Definition}\label{}\rm
Let $F=\{f_1,\ldots,f_q\}\subset R=k[{\bf x}]=k[x_1,\ldots,x_n]$ be forms of fixed
degree $d\geq 2$, with $q\geq n$. Consider a presentation of the Rees algebra ${\cal
R}_R(I)\simeq k[{\bf x}, {\bf y}]/{\cal J}$ where ${\bf y}=\{y_1,\ldots, y_q\}$ and
${\cal J}$ is a bihomogeneous ideal. We will say that $I=(F)$ is {\it of residual
linear type\/} if ${\cal J}$ is generated in bidegrees $(*,1)$ and $(0,*)$, where
$*$ denotes an arbitrary integer $\geq 1$.
\end{Definition}

\bigskip

Ideals of residual linear type are called ideals of {\it fiber type\/} in
\cite{polymatroid1}, it is shown there that polymatroidal ideals (see
Section~\ref{classified}) 
are
of fiber type. Thus, $I=(F)$ is of residual linear type if its relations are generated by the
relations that define the symmetric algebra ${\cal S}_R(I)$ and the polynomial
relations of $I$ with coefficients in the base field $k$. A conjecture -- perhaps
only a question -- regarding these ideals can be phrased as follows.

\begin{Conjecture}\label{res_lin_type}
Let $F$ be a finite set of $q\geq n$ monomials of the same degree $d\geq 2$ such
that the ideal $(F)\subset k[{\bf x}]$ is of residual linear type. Then the
following conditions are equivalent:
\begin{enumerate}
\item[{\rm (i)}] Both the log-matrix and the linear syzygy
matrix of $F$ have maximal rank.
\item[{\rm (ii)}] The extension $k[F]\subset k[{\bf x}_d]$ is
birational.
\end{enumerate}
\end{Conjecture}

A comment on the reasonableness of the conjecture. The implication (i)
$\Rightarrow$ (ii) is just Theorem~\ref{boa-viagem}.

The reverse implication (ii) $\Rightarrow$ (i) follows from the principle of linear
obstruction \cite[Proposition 3.5]{Si2} (see also \cite{CRS}) in the case of an
ideal of linear type (necessarily, $q=n$). In order to suitably extend to ideals of
residual linear type, one could in principle use the main criterion of \cite{Si2}
and the terminology thereof. Let $\phi_1$ denote the linear syzygy matrix of $F$.
Thus, the {\it weak Jacobian matrix\/} $\psi$ of $F$ (\cite[Definition 2.2]{Si2})
can be thought of as the ${\bf y}$-Jacobian matrix of the quadrics in $k[{\bf y}]$
obtained by replacing every product $x_iy_k$ in ${\bf y}\cdot \phi_1$ by $y_iy_k$,
if $1\leq i < k\leq n$, and by $(1/2)y_k^2$ if $1\leq i=k\leq n$ (thus, we need
char$(k)\neq 2$). In the case $q=n$, an easy strong duality works here to yield
that $\psi^t$ and the Jacobian dual of $\psi^t$ define the same cokernel, hence
have the same rank. But the Jacobian dual of $\psi^t$ is $\phi_1$, hence
$\rk(\psi)= \rk(\phi_1)=n-1$. For ideals of residual linear type, one needs an
analogue that says $\rk_S(\psi)=n-1 \Rightarrow \rk _{k[{\bf x}]}\phi_1=q-1$, where
$S=k[{\bf y}]/P\simeq k[F]$, with $P$ a (prime) toric ideal. Since $F$ is monomial,
a sufficiently elaborated application of Corollary~\ref{fitting_of_monomials} shows
that the minors of $\psi$ are monomials. Since $P$ is toric, then
$\rk_S(\psi)=\rk_{k[{\bf y}]}(\psi)$. Therefore, we are reduced to show that
$\rk_{k[{\bf y}]}\psi=n-1 \Rightarrow \rk _{k[{\bf x}]}\phi_1=q-1$. It is this the
missing argument, for which one may have to bring in the other underlying facts of
birationality --  e.g., the log-matrix of $F$ has maximal rank and, moreover,
coker$_S(\psi^t)$ is torsion free as $S$-module (the latter issues from the
criterion in \cite{Si2}).

\bigskip

Henceforth we assume that $\deg(x^{v_i})=2$ for all $i$. It is convenient to
interpret a set of monomials of degree two in terms of graphs, possibly with loops.
Thus, consider the graph $\widetilde{\cal G}$ on the vertex set
$X=\{x_1,\ldots,x_n\}$ whose set of edges and loops correspond bijectively to the
pairs $\{x_i,x_j\}$ such that $x_ix_j\in F$ (possibly $i=j$). Denote by ${\cal G}$
the underlying simple graph obtained by omitting all loops. Notice that, in our
situation, the log-matrix $A$ of $F$ is the incidence matrix of $\widetilde{\cal
G}$ and the monomial subring $k[F]$ is the edge subring $k[\widetilde{\cal G}]$ of
the graph $\widetilde{\cal G}$.

\medskip

One basic result for cohesive sets of monomials in degree $d=2$ reads as follows.

\begin{Proposition}\label{cohesive_vs_linear}
If $F=\{x^{v_1},\ldots,x^{v_q}\}$ is a set of forms of degree $2$ with no
non-trivial common factor. Then $\rk\,{\cal L}{\cal S}(F)=q-1$ if and only if $F$
is cohesive.
\end{Proposition}
\demo One implication follows immediately from Lemma~\ref{nonsense}. For the
reverse implication, assume that $F$ is cohesive. Then the corresponding graph
$\widetilde{\cal G}$ as above is connected, hence the underlying simple graph $\cal
G$ has a spanning tree $\cal T$. Being a tree, $\cal T$ has $n-1$ edges. The
required result is easily verified in this case by induction on the number $n$ of
vertices: consider the subtree ${\cal T}\setminus x_n$ obtained by removing a
vertex of degree one and the corresponding edge, say, $x_ix_n$. By the inductive
assumption, $\rk\,{\cal L}{\cal S}({\cal T}\setminus x_n)=n-3$, so let $L$ denote
an $(n-2)\times (n-3)$ submatrix thereof of rank $n-3$. If $x_ix_j$ is any edge of
${\cal T}\setminus x_n$ then, by restoring the removed vertex and edge, yields a
linear syzygy of ${\cal T}$ involving edges $x_ix_j$ and $x_ix_n$ and a submatrix
of ${\cal L}{\cal S}({\cal T})$ formed by bordering $L$ with the corresponding
column syzygy and a last rows of zeros. It is clear that this $(n-1)\times (n-2)$
has rank $n-2$.

This takes care of the spanning tree $\cal T$. Next, one successively restores
edges and loops on to $\cal H$ in order to recover the whole $\widetilde{\cal G}$,
this time with no new vertices. By a similar token, adding one such edge or loop at
a time to the connected subgraph $\cal H$, will increase by one the rank of the new
submatrix of ${\cal L}{\cal S}(F)$ formed by bordering as before the previous one
with the column corresponding to the added edge or loop. \QED

\medskip

Before we set ourselves to state the main result of this section, the following
observation seems pertinent. Quite generally, as used in the proof of
Lemma~\ref{nonsense} and easily shown, the field of fractions of the $d$-Veronese
algebra is generated by the fractions $x_2/x_1,x_3/x_1,\ldots,x_n/x_1$ and the pure
power $x_1^d$. Thus, a simple $\underline{\rm necessary}$ condition in order that
$k[F]\subset k[{\bf x}_d]$ be birational is that $x_1^d$ be expressed as a fraction
whose terms are products of the monomials in $F$. Now, in particular, if all these
monomials are squarefree then a reasonable {\it tour de force\/} may be needed in
order to accomplish it. Thus, e.g., for $d=2$ it is not difficult to guess that the
corresponding simple graph must have a cycle of odd length. At the other end of the
spectrum it is possible, by such elementary considerations, to guess
$\underline{\rm sufficient}$ conditions under which one has enough fractions
$x_i/x_1$ out of the monomials in $F$.

We chose to follow a more conceptual thread.

The next result generalizes \cite[Corollary~3.2]{SiVi} and gives a complete answer
for monomial birationality in degree two.

\begin{Theorem}\label{aronvila05}
Let $F\subset R$ be a finite set of monomials of degree two having no non-trivial
common factor and let ${\cal G}\subset \widetilde{\cal G}$ denote the corresponding
graphs as above. Let $A$ denote the incidence matrix of $\widetilde{\cal G}$. The
following conditions are equivalent:
\begin{enumerate}
\item[{\rm (i)}] $F$ is cohesive and $A$ has maximal rank.
\item[{\rm (ii)}] The extension $k[F]\subset k[{\bf x}_2]$ is birational.
\item[{\rm (iii)}] ${\cal G}$ is connected and, moreover, either it is non
bipartite or else it is bipartite and $\widetilde{\cal G}\setminus
 {\cal G}\neq\emptyset$.
\end{enumerate}
\end{Theorem}
 \demo We first show that (i) and (ii) are equivalent.
Clearly, (ii) implies that $\rk(A)=\dim k[F]=n$ and cohesiveness follows from
Proposition~\ref{cohesive_vs_linear}. The converse is a consequence of
Theorem~\ref{boa-viagem} and Proposition~\ref{cohesive_vs_linear}.

We next show that (ii) and (iii) are equivalent.

First, (iii) $\Rightarrow$ (ii).

  Since ${\cal G}$ is connected, there is a
spanning tree $T$ of ${\cal G}$ containing all the vertices of ${\cal G}$, see
\cite{Har}.

If ${\cal G}$ is a bipartite graph and $x_n$ is a loop of $\widetilde{\cal G}$. We may
then regard $T$ as a tree with a loop at $x_n$. Notice that $T$ has
 exactly $n-1$ simple edges plus a loop. The
incidence matrix $B$ of $T$ has order $n$, is non singular, and we may assume that
the last column of $B$ is the transpose of $(0,0,\ldots,0,2)$. Consider the matrix
$B'$ obtained from $B$ by removing the last column. The matrix $B'$ is totally
unimodular because it is the incidence matrix of a simple bipartite graph
\cite[p.~273]{Schr}. Therefore $\det(B)=\pm 2$ and ${\rm rank}(A)=n$. From
Lemma~\ref{principle2} we obtain that $k[T]\subset k[{\bf x}_2]$ is birational,
hence $k[G]\subset k[{\bf x}_2]$ is birational as well.

Now, let ${\cal G}$ be a non bipartite graph. Then ${\rm rank}(A)=n$. Since ${\cal
G}$ has a spanning tree and $G$ has at least one odd cycle (\cite[pp. 37-39 and
p.~42]{Har}), then ${\cal G}$ admits a connected simple subgraph ${\cal H}$ with
$n$ vertices and $n$ edges with a unique cycle of odd length. By
\cite[Corollary~3.2]{SiVi} the extension $k[{\cal H}]\subset k[{\bf x}_2]$ is
birational, hence  so is $k[{\cal G}]\subset k[{\bf x}_2]$.

Finally, we show the implication (ii) $\Rightarrow$ (iii).

By Proposition~\ref{cohesive_vs_linear}, $F$ must be cohesive, i.e., ${\cal G}$ is
connected. We have already seen that $\rk\,A=\dim k[F]=n$. Suppose that ${\cal G}$
is bipartite. Then the log-matrix of ${\cal G}$ has rank $n-1$, hence
$\widetilde{\cal G}$ has at least one loop. \qed

\begin{Example}\rm
A geometer would summarize the result of Theorem~\ref{aronvila05} by saying that
any cohesive coordinate projection of the $2$-Veronesean that preserves dimension
is birational onto the image. This is clearly false if the projection is a non
coordinate cohesive projection, e.g., if $F$ is a set of $2$-forms forming a
cohesive regular sequence (the simplest example with $n=2$ would be
$F=x_1x_2,x_1^2-x_2^2$). At the other end, for $d>2$, a cohesive coordinate
projection of the $d$-Veronesean preserving dimension can fail to be birational for
the simple reason that it may be composed with a non-cohesive set. The simplest
example of this phenomenon is $F=\{x_1^4, x_1^2x_2^2, x_2^4\}$. Here, $k[F]\subset
k[(x_1,x_2)_4]$ is not birational, but its ``reparametrization'' $F'=\{y_1^2,
y_1y_2, y_2^2\}$ is the $2$-Veronesean. For $n>2$, one of the simplest examples is
$F=\{x_1^3,x_1^2x_2, x_2x_3^2\}$, which is cohesive of maximal rank,
non-reparametrizable and non-birational: the ideal $(F)\subset k[x_1,x_2,x_3]$ is
of linear type, but the linear syzygy matrix is of rank $1$, hence falls below the
needed value $2$ (of course, this apparatus in such a simple example is worthless
since one immediately sees that $x_3^3$ does not belong to the field of fractions
of $k[x_1^3,x_1^2x_2, x_2x_3^2]$).
\end{Example}

\begin{Corollary}\label{contraction_of_bipartite} Let ${\cal G}$ be a 
connected simple
bipartite graph.  Assume that
$x_{n-1}x_n$ is an edge of an
even cycle of ${\cal G}$. Then $k[\widetilde{{\cal G}
\setminus x_n}]\subset k[({\bf
x}\setminus {x_n})_2]$ is a
birational extension, where $\widetilde{{\cal G}\setminus x_n}$ is the graph on
the vertices $X\setminus
{x_n}=\{x_1,\ldots,x_{n-1}\}$ obtained by contracting the
edge $x_{n-1}x_n$ to a loop
around the vertex $x_{n-1}$.
\end{Corollary}

\demo  By the contracting-looping transformation, the
resulting graph $\widetilde{{\cal G}\setminus x_n}$
acquires an odd cycle. Therefore, the simple subgraph induced by
$\widetilde{{\cal G}\setminus x_n}$ is
non-bipartite and the assertion
follows from Theorem~\ref{aronvila05}.\QED

\begin{Remark}\rm
(a) The fact that a connected graph on $n$ vertices having exactly $n$ edges and a
unique cycle of odd length induces a birational (Cremona) map had been guessed in
\cite[Conjecture 2.8]{RuSi} and proved in \cite[Corollary 3.3]{SiVi} in a
characteristic-free way. In characteristic zero, the more general context envisaged
in \cite{CRS} includes this result.

(b) If $q=n$, 
Corollary~\ref{contraction_of_bipartite} has a pretty geometric 
interpretation.
The given ring
extension $k[{\cal G}]\subset k[{\bf x}_2]$ (${\cal G}$ bipartite) 
translates into
a rational map
$${\cal F}\colon\pp^{n-1}\dasharrow \pp^{n-1}$$
whose image is Proj$(k[{\cal G}])$, after normalizing the grading of $k[{\cal G}]$.
The induced ring extension $k[\widetilde{{\cal G}\setminus x_n}]\subset k[({\bf
x}\setminus {x_n})_2]$ corresponds to the restriction of ${\cal F}$ to the
hyperplane $L$ defined by $x_{n-1}-x_n=0$ and its image can be identified with the
image of ${\cal F}$ (actually, the algebras $k[{\cal G}]$ and
$k[\widetilde{{\cal G}\setminus x_n}]$ are isomorphic 
as graded $k$-algebras by the
contracting isomorphism $k[{\bf x}]/L\simeq k[{\bf x}\setminus x_n]$ sending
$x_i\mapsto x_i$ for $1\leq i\leq n-1$ and $x_n\mapsto x_{n-1}$). 
Thus, ${\cal F}$ restricts to a birational map of $L\simeq \pp^{n-2}$ onto im$({\cal F})$.
\end{Remark}

\section{Hall of examples}\label{classified}

\subsection{Monomial Cremona transformations}

Among monomial birational maps, the Cremona ones form a well-known distinguished
class. A Cremona map is a birational map of $\pp^{n-1}$ onto itself. A recent
surprising result (\cite{PanEtAl}) showed that the monomial Cremona
transformations of $\pp^{n-1}$ is generated by the ones of degree $2$ and by the
projective linear group, thus partially extending the classical result of M.
Noether to higher dimension. The question as to which are the ``standard ones'' in
dimension $\geq 3$, if any at all, remains open as far as we know.

\subsubsection{Monomial Cremona transformations of degree $2$}

We add a tiny contribution towards further understanding the structure of such
maps. The next result extends a bit \cite[Corollary 3.3]{SiVi} and likewise
clarifies the algebraic/combinatorial background of the involved Cremona maps.

\begin{Proposition}\label{cremonadegree2}
Let $F\subset k[x_1,\ldots, x_n]$ be a cohesive finite set of monomials of degree
two having no non-trivial common factor and let ${\cal G}\subset \widetilde{\cal
G}$ denote the corresponding graphs as above. Let $A$ denote the $n\times n$
incidence matrix of $\widetilde{\cal G}$. The following conditions are equivalent:
\begin{enumerate}
\item[{\rm (i)}] $\det A\neq 0$
\item[{\rm (ii)}] $F$ defines a Cremona
transformation of $\pp^{n-1}$
\item[{\rm (iii)}] Either
\begin{itemize}
\item[{\rm (a)}] $\widetilde{\cal G}=
 {\cal G}$ {\rm (}i.e., no loops{\rm )}, ${\cal G}$ has a unique cycle
 and this cycle has odd length;
\end{itemize}
\vskip-8pt
 or else
\begin{itemize}
\item[{\rm (b)}] $\widetilde{\cal G}$ is a tree with exactly one loop.
\end{itemize}
\item[{\rm (iv)}] The ideal $(F)\subset k[x_1,\ldots,x_n]$ is of linear type.
\end{enumerate}
\end{Proposition}
\demo The equivalence of (i) through (iii) follows immediately from
Theorem~\ref{aronvila05}, by noticing that if the underlying simple graph ${\cal
G}$ is bipartite and $\widetilde{\cal G}$ has exactly $n$ edges and loops, then the
latter has to be a tree with exactly one loop. To see that the first three
conditions are also equivalent to (iv), notice that (iv) implies (i) since the
generators of an ideal of linear type are analytically independent, hence
algebraically independent as they are forms of the same degree. Now, when
$\widetilde{\cal G}= {\cal G}$, the implication (iii)(a) $\Rightarrow$ (iv) is part
of \cite[Corollary 3.3]{SiVi} but has really been noticed way before in
\cite[Corollary 3.2]{Villa2} (see also \cite[Corollary 8.2.4]{VillaBook}). Thus, it
remains to see that (iii)(b) $\Rightarrow$ (iv) in the case where $\widetilde{\cal
G}$ effectively has loops. This follows from Lemma~\ref{apr7-05} below using
induction and noticing that an edge with a loop is clearly of linear type. \qed

\begin{Lemma}\label{apr7-05} Let $F=\{f_1,\ldots,f_q\}\subset R$ be a
set of monomials of degree two and let $f_{q+1}=x_ix_{n+1}$ be a monomial in
$R'=R[x_{n+1}]$, where $x_{n+1}$ is a new variable and $1\leq i\leq n$. If $I=(F)$
is of linear type, then $I'=(I,f_{q+1})$ is of linear type.
\end{Lemma}

\demo  Let $R'[I't]$ be the Rees algebra of $I'$ over the extended polynomial ring
$R'$. Let $J'$ denote the presentation ideal of $R'[I't]$, i.e., the kernel of the
graded epimorphism:
$$
\varphi\colon\, B'=R'[t_1,\ldots,t_{q+1}]\longrightarrow R'[I't]\longrightarrow 0\
\ \ \ (t_i\longmapsto f_it).
$$
We may assume that $J'$ extends the presentation ideal $J$ of the Rees algebra
$R[It]$ over $R$ via the natural inclusion $R[It]\subset R'[I't]$. We know that
$J'=\oplus_{s\geq 1}J'_s$ is a graded ideal in the standard $\zz$-grading of
$R'[I't]$ with $(R'[I't])_0=R'$. To show that $I'$ is of linear type we have to
show that $J'_s\subset B'J'_1$ for all $s\geq 1$. We proceed by induction on $s$,
the result being vacuous for $s=1$. Thus, assume $s\geq 2$. Since $J'$ is a toric
ideal, it is generated by binomials. Therefore, by the inductive hypothesis, it
suffices to show that any binomial in $J'$ belongs to $B'J'_{s-1}$. Let
$$
h=x^\alpha t_{i_1}^{a_1}\cdots t_{i_k}^{a_k}-x^\beta
t_{j_1}^{b_1}\cdots t_{j_r}^{b_r}
$$
be a binomial in $J'_s$, where ${i_1},\ldots,i_k,j_1,\ldots,j_r$ are distinct
integers between $1$ and $n+1$, $a_i>0,b_i>0$ for all $i,j$ and $a_1+\cdots +a_k=
b_1+\cdots + b_r=s$. We may assume that $f_{i_k}=f_{q+1}=x_ix_{n+1}$, otherwise
$h\in BJ_1\subset B'J'_1$ because $I$ is of linear type. From the equality
$$
x^\alpha f_{i_1}^{a_1}\cdots f_{i_k}^{a_k}=x^\beta
f_{j_1}^{b_1}\cdots f_{j_r}^{b_r}
$$
follows that $x_{n+1}$ divides $x^\beta$, since no $f_j$ on the right side of this
equality involves the variable $x_{n+1}$. Thus there is a relation
\begin{equation}\label{relation}
x^\alpha f_{i_1}^{a_1}\cdots f_{i_{k-1}}^{a_{k-1}}f_{q+1}^{a_{q+1}-1}=x^\delta
f_{j_1}^{c_1}\cdots f_{j_r}^{c_r}
\end{equation}
where one of the $c_i$'s may be zero and $c_1+\cdots+c_r=s-1\geq 1$. We may assume
that $c_1>0$ because not all $c_i$'s are zero. Consider the equality
\begin{equation}\label{combination}
h= t_{i_{q+1}}F_1+t_{j_1}F_2,
\end{equation}
where $F_1=x^\alpha t_{i_1}^{a_1}\cdots t_{i_{k-1}}^{a_{k-1}}
t_{i_{q+1}}^{a_{q+1}-1}-x^\delta t_{j_1}^{c_1}\cdots t_{j_r}^{c_r}$ and
$F_2=x^{\delta}t_{i_{q+1}} t_{j_1}^{c_1-1}t_{j_2}^{c_2}\cdots t_{j_r}^{c_r}-x^\beta
t_{j_1}^{b_1-1}t_{j_2}^{b_2}\cdots t_{j_r}^{b_r}$. Since $F_1\in J'$ because of
(\ref{relation}), then $t_{j_1}F_2\in J'$, hence $F_2\in J'$ as $J'$ is prime.
Therefore, (\ref{combination}) expresses $h$ as an element of $B'J'_{s-1}$, as
required. \QED

\begin{Example}\rm
$\{x_1x_2,x_1x_3,x_2x_3\}$ and $\{x_1x_2,x_1x_3,x_3^2\}$ are examples of each of
the subcases (a) and (b) in Proposition~\ref{cremonadegree2}. They respectively
define the standard Cremona plane maps with $3$ distinct base points and with $2$
base points and one infinitely near point. The third type of standard Cremona map
is a double structure on one single point, hence is not monomial.
\end{Example}

\subsubsection{Squarefree monomial Cremona transformations}

Let $F\subset k[x_1,\ldots,x_n]$ be a set of $n$ squarefree monomials of degree $d$
and let $A$ denote the log-matrix of these monomials.

\medskip

For convenience, a set $F$ of monomials with no common factor defining a Cremona
transformation will be said to be a {\it Cremona set\/}. Since $F$ has no common
factor, the corresponding Cremona map determines $F$ uniquely. Likewise, we will
call the {\it inverse Cremona set\/} the set of monomials that define the inverse
map. We say that two squarefree monomial 
Cremona sets are {\it permutable\/} -- to mean ``equivalent''
in the lack of better terminology -- if they coincide up to a permutation of the
source and of the target variables. This is supposedly the equivalent of saying
that the two squarefree monomial Cremona maps are geometrically one and the same.

Obviously, for a given pair $n,d$, where $n$ is the number of variables and $d$ is
the degree of the monomials, there are a finite number of mutually non-permutable
Cremona sets with these values. Classifying means finding this complete list.

Classifying squarefree Cremona sets looks within grasp since necessarily $d\leq
n-1$. Up to permutability, the only Cremona transformation of degree $d=n-1$ in $n$
variables whose terms are squarefree monomials is the analogue of the classical
Steiner plane inversion, given by $F=\{x_1\cdots x_{n-1},\, x_1\cdots
x_{n-2}x_n,\,\ldots,\, x_2\cdots x_n\}$. For degrees $d\leq n-2$, the
classification becomes more involved. Our purpose in this part is to convey the
impact of combinatorics on birationality by examining some scattered examples for
low values of $n$ and degrees $d\leq n-2$. The case where $d=2$ was completely
covered by case (a) of Proposition~\ref{cremonadegree2}.

\medskip

Recall the following notion of combinatorial nature. If $F$ is a set of
monomials of the same degree with log-matrix $A=(a_{ij})$, its {\it dual
complement\/} is the set $\widehat{F}$ of monomials whose log-matrix is
$\widehat{A}=(1-a_{ij})$.
 The following basic principle guides us into further simplification.

\begin{Proposition}{\rm (Duality Principle)}\label{dualcomplement}
Let $F$ be a set of monomials in $n$ variables, of the same degree
$d$, with no common factor. Then $F$ is a Cremona set if and only if $\widehat{F}$ is a
Cremona set.
\end{Proposition}
\demo There is a known equality that works for all $n$ and $d$ (see for
instance \cite{escobar}):
$(n-d)\det(A)=(-1)^{n-1}d\det(\widehat{A})$. A simple proof
of this equality consists in adding the rows of $A$ to get $\det(A)=d\det(A')$, where:
$$A'=\left[
\begin{array}{ccc}
a_{1,1}&\dots & a_{1,n}\\
\vdots&  &\vdots\\
a_{n-1,1}&\dots&a_{n-1,n}\\
 1&\dots & 1
 \end{array}\right].$$
Similarly adding the rows of $\widehat{A}$ we get $\det(\widehat{A})=(n-d)\det(\widehat{A}')$,
where:

 $$\widehat{A}'=\left[
\begin{array}{ccc}
b_{1,1}&\dots & b_{1,n}\\
\vdots&\vdots&\vdots\\
b_{n-1,1}&\dots&b_{n-1,n}\\
 1&\dots & 1
 \end{array}\right].$$

Then $\widehat{A}'$ is obtained from $A'$ by subtracting the row
$\mathbf{1}=(1,\ldots,1)$ from each of the first $n-1$ rows of $A'$ and making a
change of sign at each step. Thus the determinants of $A'$ and $\widehat{A}'$ differ by
(at most) a sign given by $(-1)^{n-1}$. Thus
$$
\det(A)=d\det(A')=(-1)^{n-1}d\det(\widehat{A}')=(-1)^{n-1}\frac{d}{n-d}\det(\widehat{A}),
$$
which yields the required formula.

Now, if $F$ is Cremona then $|\det(A)|=|d|$ by DPB. Taking absolute
values, this formula yields $|\det(\widehat{A})|=|n-d|$. Thus $\widehat{F}$ is
a Cremona set by DPB. The reverse implication
is obtained by a symmetrical argument.\QED

\begin{Proposition}\label{5vars_deg3}
Up to permutation of the variables, the complete list of distinct squarefree
Cremona sets of degree $3$ in $5$ variables is as follows:
\begin{itemize}
\item $F=\{x_3x_4x_5,\,x_1x_4x_5,\, x_1x_2x_5, \,x_1x_2x_3,\,x_2x_3x_4\}$
{\rm(}${\cal DB}${\rm)}
\item $F=\{x_3x_4x_5,\,x_1x_4x_5, \,x_1x_2x_5, \,x_1x_3x_5,\,x_1x_2x_4\}$
{\rm($p\,$-involutive)}
\item $F=\{x_3x_4x_5,\,x_1x_4x_5, \,x_1x_2x_5, \,x_1x_3x_5,\,x_1x_3x_4\}$
{\rm($p\,$-involutive)}
\item $F=\{x_3x_4x_5,\,x_1x_4x_5, \,x_1x_2x_5,\, x_2x_4x_5,\,x_1x_2x_3\}$
{\rm(apocryphal)}
\end{itemize}
\end{Proposition}
\demo According to Proposition~\ref{dualcomplement}, the required complete list is
the complete list of the dual-complements. The latter is the list of all squarefree
degree $2$ Cremona sets obtained from Proposition~\ref{cremonadegree2}(a). Their
corresponding graphs are shown below:

\vspace{1.3cm}

$
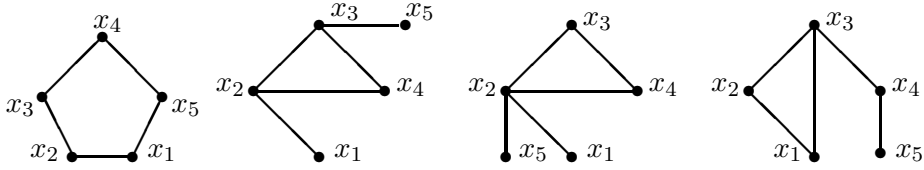
\begin{array}{ccccccc}
\setlength{\unitlength}{.040cm}
\thicklines
\begin{picture}(40,40)
\put(50,0){\circle*{3.5}}
\put(55,0){$x_1$}
\put(50,0){\line(1,2){10}}
\put(50,0){\line(-1,0){20}}
\put(30,0){\circle*{3.5}}
\put(16,0){$x_2$}
\put(30,){\line(-1,2){10}}
\put(20,20){\circle*{3.5}}
\put(8,15){$x_3$}
\put(20,20){\line(1,1){20}}
\put(60,20){\circle*{3.5}}
\put(63,15){$x_5$}
\put(60,20){\line(-1,1){20}}
\put(40,40){\circle*{3.5}}
\put(37,43){$x_4$}
\end{picture}&\ \ \ \ \ \ \ \ \ \ \ \ \ \ \ \ \  &
\setlength{\unitlength}{.025cm}
\thicklines
\begin{picture}(60,40)
\put(0,0){\circle*{5.5}}
\put(8,0){$x_1$}
\put(-35,35){\circle*{5.5}}
\put(-35,35){\line(1,0){70}}
\put(-55,35){$x_2$}
\put(0,70){\circle*{5.5}}
\put(6,75){$x_3$}
\put(35,35){\circle*{5.5}}
\put(41,35){$x_4$}
\put(-35,35){\line(1,1){35}}
\put(0,70){\line(1,-1){35}}
\put(2,70){\line(1,0){40}}
\put(46,70){\circle*{5.5}}
\put(46,75){$x_5$}
\put(-35,35){\line(1,-1){35}}
\end{picture}
 &\ \ \ \ \ \ \ \ \
&
\setlength{\unitlength}{.025cm}
\thicklines
\begin{picture}(60,40)
\put(0,0){\circle*{5.5}}
\put(8,0){$x_1$}
\put(-28,0){$x_5$}
\put(-35,35){\circle*{5.5}}
\put(-35,35){\line(1,0){70}}
\put(-55,35){$x_2$}
\put(0,70){\circle*{5.5}}
\put(6,70){$x_3$}
\put(35,35){\circle*{5.5}}
\put(41,35){$x_4$}
\put(-35,35){\line(1,1){35}}
\put(0,70){\line(1,-1){35}}
\put(-35,35){\line(1,-1){35}}
\put(-35,35){\line(0,-1){35}}
\put(-35,0){\circle*{5.5}}
\end{picture}
&
\ \ \ \ \ \ \ \
&
\setlength{\unitlength}{.025cm}
\thicklines
\begin{picture}(60,40)
\put(0,0){\circle*{5.5}}
\put(-22,0){$x_1$}
\put(43,0){$x_5$}
\put(-35,35){\circle*{5.5}}
\put(35,35){\line(0,-1){30}}
\put(35,2){\circle*{5.5}}
\put(-55,35){$x_2$}
\put(0,70){\circle*{5.5}}
\put(6,70){$x_3$}
\put(35,35){\circle*{5.5}}
\put(41,35){$x_4$}
\put(0,0){\line(-1,1){35}}
\put(-35,35){\line(1,1){35}}
\put(0,0){\line(0,0){70}}
\put(0,70){\line(1,-1){35}}
\end{picture}
\end{array}
$

\vspace{1cm}


To conclude, we explain the appended terminology. A set $F$ of squarefree monomials
is called {\it $d$-doubly-stochastic\/} (short: ${\cal DB}$) if its log-matrix
$A=(a_{ij})$ is doubly-stochastic, i.e., the entries of each column sum up to $d$
(i.e., the monomials have fixed degree $d$) and so do the entries of each row
(i.e., no variable is privileged or, the ``incidence'' degrees of the variables is
also $d$). A Cremona set is called {\it $p$-involutive\/} if it coincides with its
inverse set up to permutability. Finally, a Cremona set is called {\it
apocryphal\/} if its inverse set has at least one non-squarefree monomial. \qed

\medskip

It may be easier to classify ${\cal DB}$ Cremona sets. For instance, the following
simple result considerably reduce the possibilities.

\begin{Proposition}\label{DBobstruction} If $A=(a_{ij})$ is doubly stochastic and
$|\det(A)|=d$, then ${\rm gcd}\{n,d\}=1$.
\end{Proposition}

\demo Adding the first $n-1$ rows of $A$ to its last row and factoring out $d$ we
get:

$$\det(A)=d\det\left[
\begin{array}{cccc}
a_{11}&\cdots&a_{1n-1}&a_{1n}\\
\vdots&      &\vdots  &\vdots\\
a_{n-11}&\cdots&a_{n-1n-1}&a_{n-1n}\\
1&\cdots& 1&1
\end{array}\right].$$
Next we add the first $n-1$ columns of the matrix occurring in the right hand side
of this equality to its last column to get
$$\det(A)=d\det\left[
\begin{array}{cccc}
a_{11}&\cdots&a_{1n-1}&d\\
\vdots&      &\vdots  &\vdots\\
a_{n-11}&\cdots&a_{n-1n-1}&d\\
1&\cdots& 1 &n
\end{array}\right].$$
Hence since $\det(A)=\pm d$ we obtain that $n$ and $d$ are relatively prime. \QED

\medskip

Another useful tool is the following.

\begin{Lemma}\label{inductive_sqfree}{\rm (Inductive principle for ${\cal DB}$)}
Let $F=u_1,\ldots,u_n\subset k[x_1,\ldots,x_n]$ be a ${\cal DB}$ set of squarefree
monomials of degree $d$. Then, given a permutation $\{i_1,\ldots,i_n\}$ of the
indices such that the set $u_1/x_{i_1},\ldots,u_n/x_{i_n}$ has no repeated
monomials, then this set is a ${\cal DB}$ set of squarefree monomials of degree
$d-1$.
\end{Lemma}
\demo The proof  follows immediately from a close inspection of the corresponding
log-matrices. \qed

\medskip

Of course, the result of the lemma can be read backwards, i.e., from degree $d-1$
up to degree $d$ by multiplying by variables out of $\{x_{i_1},\ldots,x_{i_n}\}$.

\begin{Corollary}\label{noDBwith6}
For $n=6$ the only ${\cal DB}$ squarefree Cremona set of degree $d$ is the
involutive Steiner inversion, given by $$x_1x_2x_3x_4x_5,\,x_1x_2x_3x_4x_6\,,\,
x_1x_2x_3x_5x_6\,,\, x_1x_2x_4x_5x_6\,,\,x_1x_3x_4x_5x_6\,,\, x_2x_3x_4x_5x_6$$
\end{Corollary}
\demo It follows readily from Proposition~\ref{DBobstruction}. In particular the
inductive principle above does not preserve the rank of the log-matrix. \qed

\medskip

To classify the squarefree Cremona sets with $n=6$ we only need to look at degree
$d=3$, since $d=2$ follows from Proposition~\ref{cremonadegree2} and $d=4$ goes by
duality. We give some instances of Cremona sets with $n=6$ and $d=3$, with special
care for their linear syzygy behavior. These examples will hopefully give some
measure of the theoretical hardship in classifying squarefree Cremona sets for
$n\geq 6$.

\begin{Example}\rm The set $F=\{x_1x_2x_6, \,x_2x_3x_6, \,x_1x_3x_6, \,x_1x_3x_4,
\,x_1x_4x_5, \,x_3x_4x_6\}$ is a Cremona set: both the log and the linear syzygy
matrices have maximal rank. A calculation using the method of \cite{Si2} shows that
$F$ is $p$-involutive. Its dual complement $\widehat{F}$ is also a $p$-involutive
Cremona set which is not permutably equivalent to $F$.
\end{Example}

\begin{Example}\rm The set $F=\{x_1x_2x_3, \, x_2x_3x_4, \, x_3x_4x_5 ,
\, x_1x_3x_6,  \, x_2x_5x_6,\, x_4x_5x_6\}$ is a Cremona set: both the log and the
linear syzygy matrices have maximal rank. A calculation as in the previous example
shows that $F$ is apocryphal with degree $4$ inverse set $\{y_1^2y_6^2,
\,y_1y_2y_6^2,\ldots\}$ (the dots stand for squarefree monomials), a rather weird
turnout.
\end{Example}

\begin{Example}\rm The set $F=\{x_1x_2x_4, \, x_2x_3x_5, \, x_3x_4x_6 ,
\, x_1x_4x_5,\, x_1x_4x_6,\, x_2x_5x_6\}$ has log matrix of maximal rank, but not
so the linear syzygy matrix whose rank is $4$ (though the corresponding syzygy
submodule is $5$-generated). Nevertheless, a calculation as before shows that $F$
acquires an extra ${\bf x}$-linear relation (of higher ${\bf y}$-degree) which
suffices to derive birationality. Moreover, as it turns out, $F$ is apocryphal with
degree $5$ inverse set $\{y_2^2y_6^3, \,y_2^2y_3y_6^2,\,y_1y_3^2y_4y_5,\,
y_2y_3y_5y_6^2,\ldots\}$, an even weirder turnout.
\end{Example}

\subsection{Monomial birational maps from other combinatorial constructs}

The following class of sets of monomials was considered in \cite{HeHi}.

\begin{Definition}\rm A set $F=\{{\bf x}^{v_1},\ldots,{\bf x}^{v_q}\} $ of
monomials of degree $d$ minimally generating the ideal $(F)\subset k[{\bf x}]$ is
called {\it polymatroidal\/} if the following condition is satisfied: given any two
${\bf x}^u, {\bf x}^v\in F$, if $u_i> v_i$ for some index $i$ then there is an
index $j$ with $u_j < v_j$ such that $\frac{x_j}{x_i}\,{\bf x}^u\in F$.
\end{Definition}

If $F$ is polymatroidal or even matroidal, the dimension of $k[F]$
may be less than $n$. For instance if $k[F]$ is the edge subring of a
complete bipartite graph on $n$ vertices, then $\dim(k[F])=n-1$.

\medskip

The definition of polymatroidal set is somewhat tailored for having
enough
linear syzygies. This is
expressed in a slightly different way in \cite{ConHer}, where it has been shown
that, provided $F$ is ordered in the reverse lex order, it has {\it linear
quotients}, i.e., the ideals $({\bf x}^{v_1},\ldots,{\bf x}^{v_{i-1}}):{\bf
x}^{v_i}$ are generated by a set of variables, for every $i$. Clearly, this result
implies that the ideal $(F)$ is in fact linearly presented. Therefore, one has:

\begin{Proposition}\label{polymatroid}
Let $F\subset k[{\bf x}]$ be a set of monomials of degree $d$ minimally generating
the ideal $(F)$ and whose log-matrix is of maximal rank. If $F$ is polymatroidal
then $k[F]\subset k[{\bf x}_d]$ is birational.
\end{Proposition}
\demo
 This follows from Corollary~\ref{apr8-05} because $(F)$ is linearly presented as
 discussed above.
 \qed

\begin{Example}\label{veronesetype}\rm
Fix an integer $d$ and a sequence of integers $1\leq
s_1\leq\cdots\leq s_n\leq d$.
Let
$$F=\{x^{a_1}\cdots x^{a_n}\,|\,
a_1+\cdots+a_n=d;\
 0\leq a_i\leq s_i\, \forall\, i\}.$$
Then $F$ is a polymatroidal set of maximal rank (see \cite{DeHi}). The subalgebra
of $k[{\bf x}]$ generated by $F$ is said to be of {\it Veronese type}. It includes,
as special cases, the Veronese algebra of $k[{\bf x}]$ of order $d$ and the algebra
of squarefree products of $d$ variables. The birationality of $k[F]\subset k[{\bf
x}_d]$ follows directly from Proposition~\ref{polymatroid} or from \cite{SiVi}
using the fact that $R[Ft]$ is normal \cite{EVY}.
\end{Example}

\medskip

\begin{center}
ACKNOWLEDGMENT
\end{center}

\noindent The second author thank the Department of Mathematics of UFPe, where
this work started. The first author was partially supported by a CNPq
grant, Brazil.
The second author was partially supported by a CONACyT grant 49251-F 
and SNI, M\'exico.

{\small

\bibliographystyle{plain}

\begin{thebibliography}{10}

\bibitem{CRS}{C.
Ciliberto, F. Russo and A. Simis, Cremona maps, ideals of linear
type and linear
syzygies, in preparation.}

\bibitem{ConHer}{A. Conca and J. Herzog, Castelnuovo-Mumford regularity of
products of ideals, Collect. Math. {\bf 54} (2003), 137-152.}


\bibitem{DeHi}{E. de Negri and T. Hibi, Gorenstein algebras of Veronese type,
J. Algerba, {\bf 193} (1997) 629--639.}

\bibitem{Eisen}{D. Eisenbud, {\it Commutative Algebra with a view
toward Algebraic Geometry\/}, Graduate
Texts in  Mathematics {\bf 150}, Springer-Verlag, 1995.}

\bibitem{escobar}{C. Escobar, {\it Normal monomial subrings, unimodular
matrices and Ehrhart rings\/}, PhD thesis,
Cinvestav--IPN, 2004.}


\bibitem{EVY} C. Escobar, R. H. Villarreal and Y. Yoshino, Torsion
freeness and normality of blowup rings of monomial ideals,
in {\it Commutative algebra with a focus on geometric and homological
aspects\/}, Proceedings: Sevilla and Lisbon
 (A. Corso et al., Eds.), Lecture Notes in Pure and Appl. Math. 
{\bf 244}, Taylor \& Francis, Philadelphia, 2005, pp. 69-84.


\bibitem{godsil}{C. Godsil and G. Royle, {\it Algebraic Graph Theory \/},
Graduate Texts in  Mathematics {\bf 207}, Springer, New York, 2001.}


\bibitem{PanEtAl}{G. Gonzalez-Sprinberg and I. Pan, On the monomial
birational maps of the projective space, An. Acad. Brasil.
Ci\^{e}nc. 
{\bf 75} (2003), 129-134.}

\bibitem{Har}{F. Harary, {\it Graph Theory\/},
Addison-Wesley,
Reading, MA, 1972.}

\bibitem{HeHi}{J. Herzog and T. Hibi, Discrete polymatroids,
J. Algebraic Combin. {\bf 16} (2002), 239--268.}

\bibitem{polymatroid1} J. Herzog, T. Hibi and M. Vladoiu, Ideals of fiber type 
and polymatroids, Osaka J. Math. {\bf 42} (2005), 1--23.


\bibitem{RuSi}{F. Russo and A. Simis, On
birational maps and Jacobian matrices, Compositio Math. {\bf 126} (2001),
335--358.}

\bibitem{Schr}{A. Schrijver, {\it Theory of Linear and
Integer
Programming\/}, John Wiley \& Sons, New York, 1986.}

\bibitem{Si1}{A. Simis, On the
jacobian module associated to a graph, Proc. Amer. Math. Soc., {\bf 126} (1998),
989--997.}

\bibitem{Si2}{A. Simis, Cremona transformations and related
algebras,
J. Algebra {\bf 280} (2004), 162--179.}

\bibitem{SiVi} A. Simis and
R. H.  Villarreal, Constraints for the normality of monomial subrings and
birationality, Proc. Amer. Math. Soc. {\bf 131} (2003), 2043--2048.

\bibitem{Villa2}{R. H. Villarreal, Rees algebras of edge ideals, Comm. Algebra
{\bf 23} (1995), 3513--3524.}

\bibitem{VillaBook}{R. H. Villarreal, {\it Monomial
Algebras\/}, Monographs and Textbooks in Pure and Applied Mathematics {\bf 238},
Marcel Dekker, New York, 2001.}

\end{thebibliography}

}

\bigskip

\noindent {\normalsize Aron Simis} \hglue 5.6cm
{\normalsize Rafael H. Villarreal*}
\vspace{-0.75mm}\\ 
{\small Departamento de
Matem\'atica} \hglue 
2.9cm {\small Departamento de Matem\'aticas}\vspace{-1mm}\\ 
{\small Universidade Federal
de Pernambuco}\hglue 1.9cm 
{\small Centro de Investigaci\'on y de Estudios Avanzados del
IPN}\vspace{-1mm}\\ 
{\small 50740-540 Recife, Pe, Brazil}
\hglue 3.2cm {\small Apartado Postal 14--740}\vspace{-1mm}\\
{\small e-mail: {\tt aron@dmat.ufpe.br}}\hglue 3.3cm 
{\small 07000 M\'exico City, D.F.}\vspace{-1mm}\\ 
\ \hglue 7.65cm {\small e-mail: {\tt vila@math.cinvestav.mx}}

\medskip

\noindent {\small Eingegangen am 18. Juli 2005}

\end{document}